\documentclass[11pt,a4paper]{article}
\usepackage{graphicx}

\usepackage{graphicx,float,color}
\usepackage{subfigure}

\usepackage[margin=90pt]{geometry}
\usepackage{amsfonts}
\usepackage{indentfirst}             
\usepackage[sf]{titlesec}            
\usepackage{titletoc}                
\usepackage{fancyhdr}                
\usepackage{type1cm}                 
\usepackage{float}  				 
\usepackage{indentfirst}             
\usepackage{makeidx}                 
\usepackage{textcomp}                
\usepackage{layouts}                 
\usepackage{bbding}                  
\usepackage[noadjust]{cite}                   
\usepackage{color,xcolor}            
\usepackage{listings}                
\usepackage{graphics}
\usepackage{graphicx}
\usepackage{epsfig}
\usepackage{amsmath}
\usepackage{amssymb}
\usepackage{amstext}
\usepackage{subfigure}
\usepackage{amsthm}
\usepackage[runin]{abstract}
\usepackage{tikz}
\usetikzlibrary{patterns}
\usepackage{arydshln}

\usepackage{tikz}
\usetikzlibrary{patterns}
\usetikzlibrary{shapes.geometric}

\numberwithin{equation}{section}

\newtheorem{definition}{Definition}[section]
\newtheorem{theorem}[definition]{Theorem}

\newtheorem{remark}{Remark}[section]

\newtheorem{prop}[definition]{Proposition}

\usepackage[titletoc]{appendix}

\newcommand{\be}{\begin{equation}}
\newcommand{\ee}{\end{equation}}

\begin{document}

\title{Fourier series and sidewise profile control of 1-$d$ waves
}


\author{Enrique Zuazua
\thanks{Chair for Dynamics, Control, Machine Learning and Numerics, Alexander von Humboldt-Professorship, Department of Mathematics, Friedrich-Alexander University, Erlangen-N\"{u}rnberg, 91058 Erlangen, Germany.\newline (enrique.zuazua@fau.de)} 
\thanks{Chair of Computational Mathematics, Fundaci\'{o}n Deusto, University of Deusto, 48007 Bilbao, Basque Country, Spain.} 
\thanks{Departamento de Matem\'{a}ticas, Universidad Aut\'{o}noma de Madrid, 28049 Madrid, Spain.}}

\date{}

\maketitle

\begin{abstract}
We discuss the sidewise control properties of 1-$d$ waves. In analogy with classical control and inverse problems for wave propagation, the problem consists on controlling the behaviour of waves on part of the boundary of the domain where they propagate, by means of control actions localised on a different subset of the boundary. In contrast with classical problems, the goal is not to control the dynamics of the waves on the interior of the domain, but rather their boundary traces. It is therefore a goal oriented controllability problem. 

We propose a duality method that reduces the problem to suitable new observability inequalities, which consist of estimating the boundary traces of waves on part of the boundary from boundary measurements done on another subset of the boundary.   These inequalities lead to novel questions that do not seem to be treatable by the classical techniques employed in the field, such as Carleman inequalities, non-harmonic Fourier series, microlocal analysis and multipliers.

We propose a genuinely  1-$d$ solution method, based on sidewise energy propagation estimates yielding a complete sharp solution.

The obtained observability results can be reinterpreted in terms of Fourier series. This leads to new non-standard questions in the context of non-hamonic Fourier series. 
\end{abstract}

\section{Introduction and problem formulation}\label{secintro}
Wave equations intervene in one way or another in most  natural phenomena and technological applications. In the tradition of applied mathematics, waves are modelled  by hyperbolic  Partial Differential Equations (PDE), to later develop the corresponding mathematical analysis of the existence, uniqueness, regularity,  and qualitative behaviour of solutions. Most often,  when facing real applications, one needs to go further, controlling the dynamics, so to ensure the correct performance of the mechanisms and dynamics under consideration. This requires further understanding of the intrinsic properties of the dynamics and its solutions. Of course, numerical analysis and scientific computing can contribute to this endeavour. But increasing difficulties are encountered  when trying to capture high-frequency wave phenomena (\cite{E2}). Furthermore, in some applications, the complete model is not known and this leads to inverse and parameter identification problems.

All this constitutes an ambitious program to which significant efforts have been devoted, developing and combining various tools from analysis such as asymptotic methods, multipliers, Carleman inequalities, Fourier analysis, microlocal analysis, etc.

Nowadays, when addressing these problems, one can also enjoy of the added tools that Machine Learning provides, with a significant and rich interface with control theory (see \cite{RZ}). But in these notes, we will remain in the classical context of Fourier and PDE analysis of the control of waves.

A first question arises when addressing all these issues: Should one analyse waves in the physical space or in the Fourier one? And the answer is, of course, aligned with the classical wave-particle duality: a complete understanding requires combining both approaches.

The convenience of systematically using  Fourier series techniques to address control problems for waves (and also for heat-like processes) was early identified by the PDE control community. The classical survey article by D. L. Russell \cite{Ru} gathers some of the early developments in that context. 

Later on, Yves Meyer greatly contributed to expanding the in-depth use of the theory of non-harmonic Fourier series for wave-like control.  The motivation  and genesis of the contributions of Yves Meyer to this field are due to a large extent to his interaction with Jacques-Louis Lions.   Lions, who in the late 80s developed a systematic analysis of the controllability properties of various wave-like models, and who rapidly observed  the need for Fourier analysis tools to disclose some of the most intricate issues. As he used to do when looking for expert advice, he approached  Yves Meyer with a question related to the control of the plate equation. Lions had observed that  plate models, being of order four in space, enjoyed some added control properties that were not fulfilled by the wave equation. In particular, being well aware of the intrinsic  infinite velocity of propagation of plate models, he observed that, in some cases, plates could be controlled in arbitrarily short time intervals and from subsets of the boundary that were significantly smaller than those needed to control the standard second order wave model. 

This motivated Meyer's interest whose contributions to the control of wave-like  models are summarised in the synthesis article by A. Cohen \cite{Cohen} (Section 4.1). In particular,    in \cite{Meyer1}, \cite{Meyer2}, and \cite{Meyer3}, Meyer  discloses some of the fine properties of solutions of the wave equation, such as the existence of levitating solutions, resonance phenomena, and sharp regularity. 

In his interaction with Lions, Meyer rapidly realised and observed the potential depth of the question related to the plate model and he transmitted it to St\'ephane Jaffard,  suggesting him to study the works of Jean-Pierre Kahane on the pseudo-periodicity of lacunary Fourier series. Later, St\'ephane Jaffard got in contact with Alain Haraux, my PhD advisor, who was also by then intensively communicating with Lions on related topics. They collaborated to produce the relevant and influential article \cite{HJ}. In this paper, they developed an extremely interesting corollary of the celebrated Beurling-Malliavin's Theorem allowing them to prove the point-wise spectral controllability of the plate model in domains where the spectrum of the underlying elliptic generator is simple, a property that is fulfilled generically with respect to the shape of the plate. The fundamental result in \cite{HJ} on non-harmonic Fourier series replaced the spectral gap condition in the classical Ingham's inequality by a much weaker condition on the high-frequency asymptotic density of the spectrum, of a much wider validity,  to guarantee spectral controllability properties, weaker than the classical control results in energy spaces, but actually the natural and only possible one for many other problems of vibrations in which Ingham's gap condition is not fulfilled.

The results and methods \cite{HJ}  could later be used in many other contexts, such as, for instance, when dealing with wave propagation on networks, \cite{Dager}. 

This passage shows the generosity and exceptional vision of Yves Meyer to promote and pursue fruitful scientific connections and collaborations, influencing generously  scientists in different areas and generations. 
The interaction between Jacques-Louis Lions and Yves Meyer was particularly rich in the interface between PDE control and Fourier analysis. There are many other examples of how influential  this dynamics became. For instance, in \cite{JTZ}, we gave a generalisation of Ingham's inequality to families of non-harmonic Fourier series that, under a weaker gap condition,  allows us to get quantity versions of spectral inequalities.

It would be impossible to summarise here the state of the art in the theory of non-harmonic Fourier series and how it has evolved to face the challenges arising in the control of vibrations. But it is worth recalling what the fundamental Ingham inequality assures.

\begin{theorem}[Ingham \cite{i}, \cite{Y}]  \label{te.i1}
Let $(\lambda_n)_{n\in\mathbb{Z}}$
be a
sequence of real numbers and $\gamma>0$ be such that
\begin{equation}\label{eq.co1i}
\lambda_{n+1}-\lambda_n\geq \gamma>0,\quad \forall n\in
\mathbb{Z}.
\end{equation}

For any real $T$ with
\begin{equation}\label{eq.co2i}
T> \pi/\gamma \end{equation} there exist  positive
constant $C_1=C_1(T,\gamma)>0$ and $C_2=C_2(T,\gamma)>0$ such that, for any finite sequence
$(a_n)_{n\in\mathbb{Z}}$,
     \begin{equation}\label{eq.des1}
C_1\sum_{n\in\mathbb{Z}}|a_n|^2\leq \int^T_{-T}\left|\sum_{n\in\mathbb{Z}} a_n
e^{i\lambda_nt}\right|^2dt \leq C_2\sum_{n\in\mathbb{Z}}|a_n|^2.
     \end{equation}
\end{theorem}

The term of non-harmonic Fourier series refers to the fact that frequencies involved in this series of complex exponentials, $\lambda_n, n\in\mathbb{Z}$, are not integers or, more generally, equi-spaced. Ingham's inequality can be viewed as generalisation of Parseval's identity to series fulfilling the gap condition. The lower bound assures that the $\ell^2$ norm of the family of frequencies can be bounded above by the $L^2$-norm of the series in the physical space provided the time-interval is long enough, namely, $T> \pi/\gamma $. This minimal time is sharp for Ingham's inequality to hold with such a degree of generality.

These series arise naturally in the study of the control of waves when expanding solutions on the basis of the spectrum of the underlying elliptic operator. And, essentially, except for the constant coefficient d'Alembert model with the simplest boundary conditions, this leads to such series where the frequencies are not equi-spaced.

In this article we discuss the novel problem of the sidewise control of 1-$d$ waves. It can be solved using duality and sidewise energy propagation arguments, as we shall see, leading to new observability estimates for 1-$d$ waves. But, surprisingly, these results do not seem to be achievable by Fourier techniques. This article is devoted to  reinterprete these observability inequalities in the Fourier context to formulate some open problems that, hopefully, will lead to interesting developments in the future. 

``Sidewise" here and in the sequel stands for the very fact that the 1-$d$ wave equation is also well-posed when swapping the role of space and time variables.  This is consistent with the beautiful space-time symmetry  of the d'Alembert equation, which establishes the identity between the curvature and acceleration of vibrating strings. Sidewise estimates is a classical tool that has been used often in the analysis of 1-$d$ waves when dealing with its propagation and oscillatory properties,  inverse and control problems (\cite{CazHar}, \cite{Symes}, \cite{Haraux}, \cite{E1}). Here we adapt it to deal with sidewise control problems.

The techniques developed in this paper are limited to the 1-$d$ setting. We refer to \cite{Sarac} for a preliminary discussion of this topic. Similar problems can be formulated also for fourth order beam models, but we lack of tools for to treat them since the corresponding equations are ill-posed in the space-like direction.
The same sidewise control problems are relevant and can be formulated in the multi-dimensional case. But their analysis requires an in-depth use of microlocal techniques, \cite{DZ}. The parabolic counterpart of  these problems was addressed in \cite{Barcena}, also in 1-$d$, using transmutation and flatness methods in Gevrey spaces. 
\section{The sidewise control problem for 1-$d$ waves}

Consider the following variable coefficient controlled 1-$d$ wave equation:  
\begin{equation}\label{e1}
\left\{~
\begin{aligned}
&\rho (x)y_{tt}-(a(x)y_{x})_{x}=0,  && 0< x< L,~ 0< t<T,\\
&y(x,0)=y_0(x), ~y_{t}(x,0)=y_1(x),&& 0< x< L,\\
&y_x(0,t)=u(t), ~ y_x(L,t)=0,                 &&0<t<T.
\end {aligned}
\right.
\end{equation} 

In \eqref{e1}, $0<T<\infty$ stands for the length of the time-horizon, $L$ is the length of the string where waves propagate, $y=y(x,t)$ is the state and $u=u(t)$ is a control that acts on the system through the extreme $x = 0$.

This is the Neumann version of the sidewise control problem for 1-$d$ waves considered in \cite{Sarac} in the context of the Dirichlet problem.

As in \cite{Sarac}, we assume that the coefficients $\rho$ and $a$ are in $BV$,  uniformly bounded above and  below by positive constants, i. e.
\begin{equation}\label{e43}
\rho, a \in BV(0,L),
\end{equation}
\begin{equation}\label{e42}
0<\rho_0 \le \rho (x) \le \rho_1,~~0<a_0 \le a (x) \le a_1 ~~ \text{a.e. in} ~~(0,L).
\end{equation}

For any given $(y_0,y_1)$ with $y_0 \in H^1(0,L)$ and $y_1 \in L^{2}(0,L)$ and any $u \in L^2(0,T)$, system \eqref{e1} admits a unique  finite energy solution
\begin{equation}\label{transpo}
y \in C([0,T];H^1(0,L)), ~~ y_{t} \in C([0,T];L^{2}(0,L)).
\end{equation}

The problems we shall consider also make sense for bounded, positive, measurable coefficients. But  the $BV$ regularity assumption on the coefficients is, roughly, the minimal one required to guarantee that waves behave according to our intuition, very much based on d'Alembert's formula and Fourier series representations, so that they effectively travel in space-time with a uniform velocity along characteristics.
In fact, for low regularity coefficients, the high-frequency wave behavior can be rather complex \cite{FZ1}, and sidewise propagation of waves can experience loss of regularity phenomena.

\begin{figure}[h] \label{fig3}
\begin{center}
\resizebox{8cm}{!}{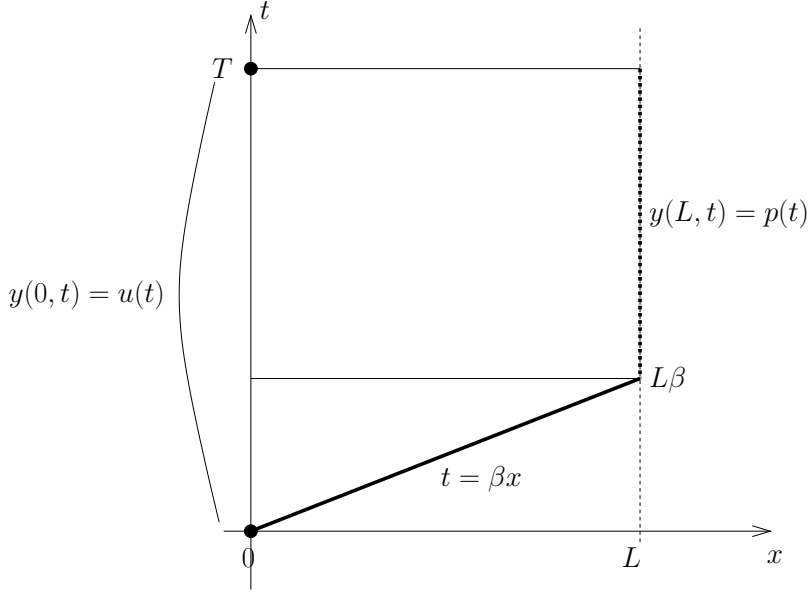}
\end{center}
\caption{The control $u(t)$ acting at $x=0$ during the time interval $[0, T]$ assures that the trace at $x=L$ reaches the target $p(t)$ within the time interval $(L\beta, T)$.}
\end{figure}

The sidewise boundary controllability property of \eqref{e1} can be formulated as follows (see Figure \ref{fig3}):

\begin{theorem}\label{controltheorem}
Let us consider system \eqref{e1} with coefficients satisfying the assumptions
\eqref{e43} and \eqref{e42}.

Let 
\begin{equation}\label{e47}T>L\beta=L \operatorname{ess}\sup\limits_{x\in[0, L]}\sqrt\frac{\rho}{a}.
\end{equation}

Then, for any $p\in H^{1}(L\beta, T)$ with $p(L\beta)=0$,  there exists a control $u\in L^2(0, T)$  such that the solution of \eqref{e1} satisfies
\begin{equation}\label{e31}
y(L, t)=p(t), \, \hbox{ for all } \, t \in (L\beta, T).
\end{equation}

\end{theorem}

\begin{remark}
To better understand this statement it is convenient to take into account the following facts:
\begin{itemize}
\item In the present formulation of the sidewise controllability property the finite velocity of propagation of waves plays an important role. 

On the one hand,   the time horizon is required to be large enough: $T>L\beta$.  This assures that waves emanating from $x=0$, where the control is active, reach the opposite extreme, $x=L$, where the tracking condition \eqref{e31} is aimed. 

On the other hand, the tracking condition is only assured in the time sub-interval $(L\beta,T)$. This sub-interval is sharp since for $t<L\beta$ the action of the control at $x=0$ has no effect on $x=L$.

\item This is a sidewise controllability result in which the action on one extreme of the string, at $x=0$, through the control $u(t)$, guarantees that the dynamics at the free end, $x=L$, is the one chosen a priori, $p=p(t)$. This is in contrast with the classical controllability problems in which the goal is to drive the vibrations, everywhere in space,  to a given configuration at the final time $t=T$ by means of suitable controls (\cite{L1}, \cite{L2}, \cite{M-E}, \cite{survey}). 

\item Without loss of generality, using the principle of additive superposition of solutions of linear systems, the problem can be reduced to the particular case $y_0(x)\equiv y_1(x)\equiv 0$.  In the sequel, our analysis will be reduced to this case.

\item 
The problem makes sense since solutions of \eqref{e1} with BV coefficients  fulfill the added boundary regularity condition 
\begin{equation}\label{boundreg}
y(L,t) \in H^{1}(0, T).
\end{equation}

\end{itemize}
\end{remark}

In \cite{Sarac} this problem was addressed in the case of Dirichlet boundary conditions. The methodology developed in that paper involved two ingredients. First, the dual equivalent sidewise observability problem was introduced, which was then solved using sidewise energy estimates. 

Note that employing the dual sidewise observability equivalent formulation in this context constitutes a natural extension of the duality methods introduced by J.-L. Lions in \cite{L1} and \cite{L2} and, in particular, the so-called Hilbert Uniqueness Method (HUM). On the other hand, in the present hyperbolic context, it is very natural to expect the same kind of properties to hold when reversing the role of space and time variables since information propagates along characteristics that are oblique in $(x, t)$. It is, therefore,  natural that methods like HUM, derived to deal with classical control problems where the goal is to control the space-like profile of the solution at a final time, might also be adapted  to the sidewise control of traces. This viewpoint led to the analysis in this paper, aimed to reformulate the problems from the non-harmonic Fourier series perspective, which leads to some interesting open problems.

Sidewise control problems are very natural in different contexts. In \cite{Gugat}, it was introduced and addressed in the context of gas flow on networks,  later extended in \cite{Li1}, \cite{Libook}, and \cite{ZLL}   to 1-$d$ quasilinear hyperbolic systems employing the method of characteristics.  In fact, sidewise control is a natural and ubiquitous concept that can be viewed as a goal-oriented control problem in which one is not aiming to fully control the trajectory but just some specific traces. In the context of structural vibrations, it may aim at preserving some parts of the structure  in a given vibrational mode. In the context of diffusion of pollutants or populations, it may refer to the need to preserve some boundary regions away from the effect of diffusion.

\section{The dual sidewise observability problem}
Let us now consider the adjoint system:
\begin{equation}\label{e45}
\left\{~
\begin{aligned}
&\rho (x)\psi_{tt}-(a(x)\psi_{x})_{x}=0, && 0< x< L,~ 0< t<T\\
&\psi(x,T)=0, ~ \psi_{t}(x,T)=0,            && 0< x<L\\
&\psi_x(0,t)=0, ~ \psi_x(L,t)=s(t),                && 0< t<T
\end {aligned}
\right.
\end{equation} \\
where the boundary data are of the form 
 \begin{equation}\label{e46}
s(t) \in H^{-1}(0, T),\,  \hbox{supp}(s) \subset (L\beta, T).
\end{equation}

This system admits a unique weak solution $\psi$ such that 
\begin{equation*}
(\psi, \frac{\partial \psi}{\partial {t}}) \in C([0,T],L^{2}(0,L) \times H^{-1}(0,L))
\end{equation*}
and
\begin{equation*}
\psi(0,.) \in L^2(0,L). 
\end{equation*}

Note that the boundary conditions in this non-standard adjoint system at the extreme $x=L$ are inhomogeneous, which is consistent with the fact that we are dealing with sidewise problems.

\begin{figure}[h] \label{fig1}
\begin{center}
\resizebox{8cm}{!}{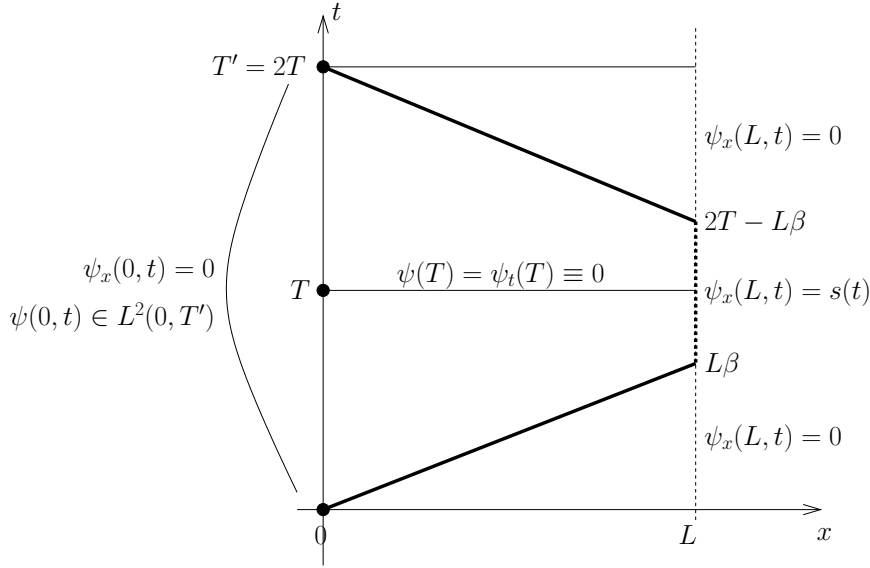}
\end{center}
\caption{Sidewise propagation of the energy of $\psi$ from $x=0$ (in the time interval $[0, T'=2T]$) towards $x=L$, to reach the time interval $[L\beta, T'-L\beta]$ at $x=L$. The solution $\psi$ is extended in an even manner with respect to time,  beyond $t=T$, so to get sharp results in terms of the regions of influence of the observed quantities at $x=0$.
}
\end{figure}

The sidewise  observability inequality that, as we shall see, this adjoint system fulfills, and which is equivalent to the sidewise controllability problem under consideration, is the following (see Figure \ref{fig1}):
\begin{prop}\label{pro41}
Assume that the coefficients $a$ and $\rho$ satisfy the assumptions
\eqref{e43} and \eqref{e42}.

 Let $T>L\beta$  as in \eqref{e47}. 

Then, there exists $C>0$ (depending on $a$, $\rho$ and $T$ but independent of $s=s(t)$),  such that the observability inequality
\begin{equation}\label{e48}
\left \Vert s(t)\right\Vert
_{H^{-1}\left(L\beta,T\right) } \leq C \left \Vert \psi(0,t)\right\Vert
_{L^{2}\left(0,T\right) } 
\end{equation}
is satisfied for every finite solution of \eqref{e45}, with $s$ as in  \eqref{e46}.
\end{prop}

\begin{remark}
It is convenient to take into account the following facts:
\begin{itemize}

\item Similar observability inequalities can also be proved in other Sobolev norms. For instance, it also holds:
\begin{equation}\label{e48bbb}
\left \Vert s(t)\right\Vert
_{L^2\left(L\beta,T\right) }\leq C \left \Vert \psi(0,t)\right\Vert
_{H^1\left(0,T\right) }.
\end{equation}

\item This result asserts that the lateral source of waves $s=s(t)$ at the extreme $x=L$  can be estimated out of the displacement $\psi(0, t)$ at the opposite extreme $x=0$.

\item The proof of this observability inequality can be easily achieved using sidewise energy estimates as in \cite{Sarac} (see also \cite{E1} and \cite{Cara}). It consists on viewing system \eqref{e45} in a sidewise manner, so that $x$ plays the role of the time variable and $t$ the space one. Classical energy estimates in this reversed sense yield the desired observability inequality.

To be more precise, consider the sidewise energy
\begin{equation*}
F(x)=\frac{1}{2}\int_{\beta x}^{T-\beta x}[\rho (x) (\psi_{t}(x,t))^2+a(x)(\psi_{x}(x,t))^2]dt, ~~~~\forall x \in [0,L].
\end{equation*}
Note that, due to the homogeneous Neumann boundary condition at $x=0$, 
$$F(0)=\frac{1}{2}\int_{0}^{T}\rho (0^+) \psi_{t}(0,t)^2 dt
$$
and that $F(L)$ yields an upper bound on the value $\psi_{x}(L,t) =s(t)$ at $x=L$.  

The needed observability inequality is therefore a consequence of an estimate of the form
$$
F(L) \ C F(0)
$$ that can be easily achieved by computing the $x$-derivative of $F$ and applying Gronwall's inequality, making use of the fact that the wave equation is fulfilled. This is in fact a classical energy estimate for the solutions of the adjoint wave equation, applied in the sidewise sense,  reversing the sense of space and time, i. e.
$$
(a(x)\psi_{x})_{x}-\rho (x)\psi_{tt}=0.
$$

We omit the details of the proof, which can be found in \cite{Sarac}, and that has been previously implemented in other 1-$d$ wave control problems as in \cite{E1} and \cite{Cara}.

Our goal in this article is to present some extensions and to analyze them from the Fourier series perspective.

\item All terms in the inequality \eqref{e48} make sense. In fact, the sidewise energy estimates in \cite{Cara} allow showing that there exists another constant $C=C(\rho,a,T)>0$ such that the finite energy solution $\psi$ of  \eqref{e45} satisfies 
\begin{equation}\label{e91}
 \left \Vert \psi(0,t)\right\Vert
_{L{2}\left(0,T\right) } \leq C \left \Vert s(t)\right\Vert
_{H^{-1}\left(0,T\right) }
\end{equation}
for all $s \in H^{-1}(0, T)$. 

\item  $BV$ is the minimal requirement on the regularity of the coefficients for this inequality to hold. Counterexamples can be built in the class of H\"older continuous ones (\cite{CZ}). In case coefficients are slightly less regular than $BV$ one may obtain weaker controllability properties with a loss of a finite  number of derivates in Sobolev-norms, \cite{FZ1}.

\end{itemize}
\end{remark}

This observability inequality \eqref{e48} is equivalent to the sidewise controllability property in Theorem \ref{controltheorem}.
In fact, out of the observability property above, one can obtain the sidewise control of minimal $L^2(0, T)$-norm by a variational principle that we present now, inspired by the so-called Hilbert Uniqueness Method (HUM) by J.-L. Lions in \cite{L1} and \cite{L2}. 

 Let us consider the continuous and strictly convex quadratic functional
 \begin{equation*}
J:H^{-1}(L\beta,T) \longrightarrow \mathbb{R}
\end{equation*}
defined as
\begin{equation}\label{e410}
J(s)=\frac{a(0^+)}{2}\int_{0}^{T}[\psi(0,t)]^2dt+a(L^-) \int_{L\beta}^{T}s(t)p(t)dt
\end{equation}
defined for  $s\in H^{-1}(L\beta,T)$ as in \eqref{e46} (i. e. with $s$ extended by zero outside $(L\beta,T)$), where $\psi$ is the corresponding solution of the adjoint system \eqref{e45}.

Note that, in \eqref{e410}, the integral $\int_{L\beta}^{T}p(t)s(t)dt$ represents the duality pairing between $ p(t) \in H^{1}_0(L\beta,T)$ and $s(t) \in H^{-1}(L\beta,T)$, i.e. if $s(t) = S'(t)$ with $S\in L^2(L\beta,T)$, then 
$$
\int_{L\beta}^{T}s(t)p(t)dt = - \int_{L\beta}^{T}S(t)p'(t)dt.
$$

The observability inequality above guarantees that the functional $J$ is also coercive. The Direct Method of the Calculus of Variations then ensures that $J$ has a unique minimizer.

The following lemma states that the minimum of the functional $J$ provides the desired sidewise control.
\begin{prop}\label{lemma1}
Suppose that  $\bar s \in H^{-1}(L\beta,T)$ is the unique minimizer of $J$ in the space $ H^{-1}(L\beta,T)$ (extended by zero outside $(L\beta,T)$). If $\bar \psi$ is the solution of the adjoint system \eqref{e45} corresponding to this boundary datum  $\bar s = \bar s(t)$, then
\begin{equation}\label{e49}
u(t)=\bar \psi(0,t) \in L^2(0, T)
\end{equation}
is the control of minimal $L^2(0, T)$-norm such that 
\begin{equation}
y(L,t)=p(t), \quad L\beta\le  t\le T,
\end{equation}
when the initial data $y_0\equiv y_1\equiv 0$.
\end{prop}
\begin{remark} 
The following comments are in order:
\begin{itemize}

\item As mentioned above, once the control is built for $y_0\equiv y_1\equiv 0$, using the linear superposition of solutions of the wave equation, the control for arbitrary initial data can be built. The functional $J$ above can also be modified so to lead directly to the control corresponding to non-trivial initial data.

\item The variational characterization of the minimal $L^2(0,T)$-norm control is also useful when dealing with numerical approximation issues. We refer to the pioneering work of R. Glowinski and J. L. Lions \cite{GL} on this topic and to \cite{E2} for further developments.
\end{itemize}
\end{remark} 

\noindent \begin{proof}
Since $J$ achieves its minimum at  $\bar s$, the Gateaux derivative of $J$ vanishes at that point. In other words, 
\begin{equation}\label{e37}
a(0^+) \int_{0}^{T}\bar \psi(0,t)\psi(0,t)dt+a(L^-)\int_{L\beta}^{T}s(t)p(t)dt=0
\end{equation}
for any other $s \in H^{-1}(L\beta,T)$ (extended by zero outside $(L\beta,T)$), where $ \psi$ stands for the solution of adjoint system \eqref{e45} with $s$ as boundary datum.

On the other hand, multiplying the equation \eqref{e1} by the solution $\psi(x,t)$ of the adjoint problem and integrating on $(0,L)\times (0,T)$, we get (recall that, without loss of generality, we have assumed that $y_0(x)\equiv y_1(x)\equiv 0$)
\begin{equation}\label{e38}
a(0^+) \int_{0}^{T}u(t)\psi(0,t)dt+a(L^-) \int_{L\beta}^{T}s(t)y(L,t)dt=0.
\end{equation}

Comparing \eqref{e37} and \eqref{e38} we get that $u(t)$ as in \eqref{e49}  is a control which leads $y(x,t)$ to $p(t)$ on $x=L$ for all $L\beta< t<T$.

Using classical arguments, it can be seen that the control obtained by this minimization method is the one of minimal $L^2(0, T)$-norm (see \cite{M-E}).

\end{proof}

\section{Sidewise plus exact control}

As we have seen, system \eqref{e1} is sidewise controllable when $T>L\beta$. On the other hand, the system is well-known to be controllable in the classical sense when $T>2L\beta$ (\cite{Cara}). It is then natural to analyse whether the system can be simultaneously controlled by the same control $u=u(t)$ in both ways, guaranteeing that the desired trace is achieved at the free end $x=L$ and the final configuration of the solution is achieved at $t=T$. As we shall see, this can be done, indeed, when $T>2L\beta$.

\begin{figure}[H] \label{fig2}
\begin{center}
\resizebox{8cm}{!}{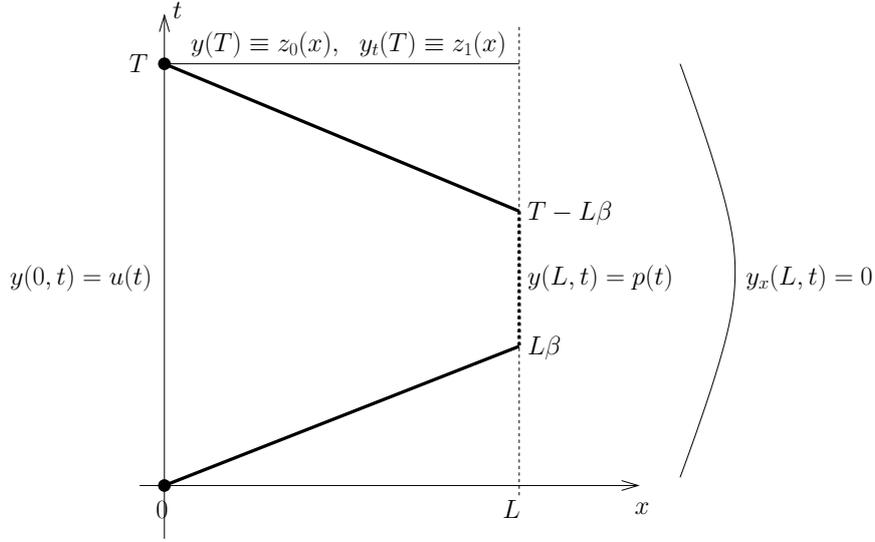}
\end{center}
\caption{The control $u(t)$ acting at $x=0$ during the time interval $[0, T]$ assures both that the trace at $x=L$ reaches the target $p(t)$ within the time interval $(L\beta, T-L\beta)$ and, simultaneously, the final conditions $(z_0, z_1)$ at $t=T$. }
\end{figure}

To make the problem more precise, in the context of system \eqref{e1}, other than the sidewise target \eqref{e31}, the  following final condition is imposed
\begin{equation}\label{final}
y(x, T) = z_0(x), y_t(x, T)=z_1(x), \quad 0<x<L.
\end{equation}
The following result holds (see Figure \ref{fig2}):
\begin{theorem}\label{esc}
Let us consider system \eqref{e1} with coefficients satisfying the assumptions
\eqref{e43} and \eqref{e42}.

Let $$T>2L\beta$$ with $\beta$ as in \eqref{e47}.

Then, for any $p\in H^{1}_0(L\beta, T-L\beta)$, $z_0\in H^1(0, L)$ and $z_1\in L^2(0, L)$  there exists a control $u\in L^2(0, T)$  such that the solution of \eqref{e1} satisfies simultaneously \eqref{final} and
\begin{equation}\label{e31b}
y(L, t)=p(t), \, \hbox{ for all } \, t \in (L\beta, T-L\beta).
\end{equation}
\end{theorem}

\begin{remark}
Several remarks are in order:
\begin{itemize}
\item Note that in Theorem \ref{esc} the time of control is twice the one in Theorem \ref{controltheorem}. This is natural since the time needed to guarantee \eqref{final} is precisely $T>2L\beta$.

\item Theorem above constitutes an improvement of classical controllability results since, other than \eqref{final}, the control also assures \eqref{e31b} to hold.

\item Note however that the sidewise control property in \eqref{e31b} is only assured to hold for $ t \in (L\beta, T-L\beta)$. As explained in the previous sections, the fact that it only holds for $t >L\beta$ is natural since this is the waiting time for the action at the controlled extreme $x=0$ to reach the other one $x=L$.  On the other hand, for the sidewise target to be reached we also need to impose the condition $p\equiv 0$ in $(T-L\beta, T)$.
\end{itemize}
\end{remark}

As in previous sections, the problem can be reduced by duality, to an observability inequality for the adjoint system that, this time, takes the form
\begin{equation}\label{e45bc}
\left\{~
\begin{aligned}
&\rho (x)\psi_{tt}-(a(x)\psi_{x})_{x}=0, && 0< x< L,~ 0< t<T\\
&\psi(x,T)=\psi_0(x), ~ \psi_{t}(x,T)=\psi_1(x),            && 0< x<L\\
&\psi_x(0,t)=0, ~ \psi_x(L,t)=s(t),                && 0< t<T
\end {aligned}
\right.
\end{equation} \\
where the boundary data are of the form 
 \begin{equation}\label{e46b}
s(t) \in H^{-1}(0, T),\,  \hbox{supp}(s) \subset (L\beta, T-L\beta),
\end{equation}
and
\begin{equation}
\psi_0(x) \in L^2(0, L), \, \psi_1(x) \in H^{-1}(0, L).
\end{equation}

Note that the main difference with respect to \eqref{e45}, is that the data of the adjoint state at $t=T$ are inhomogeneous and that $s(t)$ is assumed to vanish away from $(L\beta, T-L\beta)$.

\begin{figure}[H]\label{fig4} 
\begin{center}
\resizebox{8cm}{!}{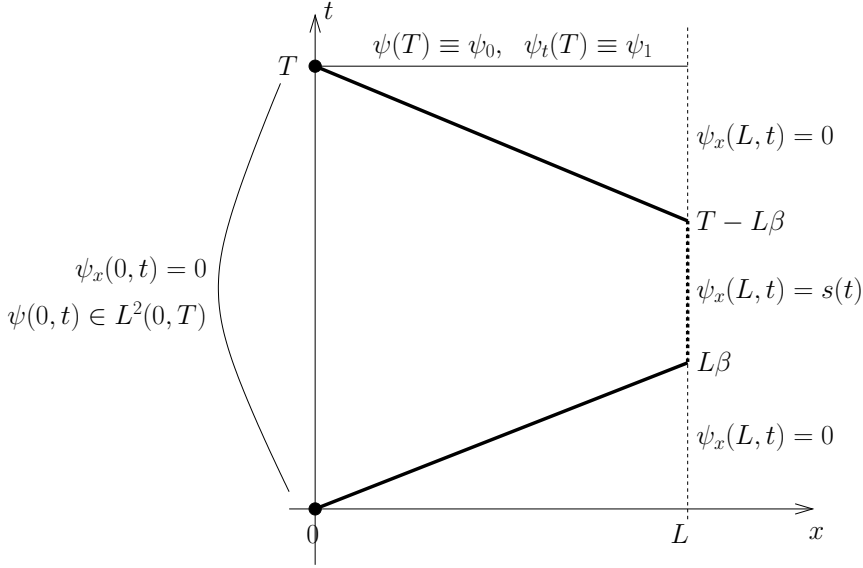}
\end{center}
\caption{Sidewise propagation of the energy of $\psi$ from $x=0$ (in the time interval $[0, T]$) towards $x=L$, to reach at $x=L$ the time interval $[L\beta, T-L\beta]$ and simultaneously recover the energy of the final data at $t=T$. }
\end{figure}

The following observability inequality holds (see Figure \ref{fig4}):
\begin{prop}\label{pro41b}
Assume that the coefficients $a$ and $\rho$ satisfy the assumptions
 \eqref{e43} and \eqref{e42}.

 Let $T>2L\beta$ ($\beta$ being given as in \eqref{e47}). 

Then, there exists $C>0$ such that the observability inequality
\begin{equation}\label{e48b'bbbb}
\left \Vert s(t)\right\Vert
_{H^{-1}\left(L\beta,T-L\beta\right) } + \left \Vert  \psi_0\right\Vert_{L^2(0, L)} + \left \Vert \psi_1 \right\Vert_{H^{-1}(0, L)}  \leq C \left \Vert \psi(0,t)\right\Vert
_{L^{2}\left(0,T\right) } 
\end{equation}
is satisfied for every finite solution of \eqref{e45bc}, with $s$ as in  \eqref{e46b}.
\end{prop}

\begin{remark} Several remarks are in order
\begin{itemize}
\item The proof of \eqref{e48b'bbbb} can be done, again, by sidewise energy estimates. But one has to proceed more carefully in two steps.
\begin{itemize}

\item Step 1. One first proceeds by sidewise energy estimates so to get the estimate
\begin{equation}\label{e48b'bbbbc}
\left \Vert s(t)\right\Vert
_{H^{-1}\left(L\beta,T-L\beta\right) } \leq C \left \Vert \psi(0,t)\right\Vert
_{L^{2}\left(0,T\right) } 
\end{equation}
This also yields an estimate on the energy of the solution  at time $t=T-L\beta$, i. e. 
\begin{equation}\label{e48b''}
 \left \Vert  \psi(\cdot, T-L\beta) \right\Vert_{L^2(0, L)} + \left \Vert \psi_t(\cdot, T-L\beta)  \right\Vert_{H^{-1}(0, L)}  \leq C \left \Vert \psi(0,t)\right\Vert
_{L^{2}\left(0,T\right) } 
\end{equation}

\item Step 2. Once this is done, because \eqref{e48b''} holds and the the boundary condition $s=s(t)$ vanishes after $t=T-L\beta$, one can simply use the energy conservation identity from $t=T-L\beta$ to $t=T$ to get \eqref{e48b'bbbb}.
\end{itemize}

\item Most of the remarks in the other sections of this paper apply in this context too, in particular in what concerns the minimal $BV$ regularity required on the coefficients.
\end{itemize}
\end{remark}

\section{Fourier series representation: An open problem}

Our goal in this section is to discuss to which extent these issues can be treated using Fourier series representations, making a link with the theory of non-harmonic Fourier series. We focus on Theorem \ref{controltheorem} and the adjoint system \eqref{e45}.

Let us introduce the spectrum of the corresponding Sturm-Liouville problem
      \begin{equation}
      \label{eqdeltalambda}
      \left\{\begin{array}{ccc}
   - (a(x) \partial_x w_k)_x=\lambda_k \rho(x) w_k & \hbox{in } (0, L)         \\
            \partial_x  w_k=0 & \quad\quad\quad\quad\hbox{ for } x=0, L,\,  k\ge 1.
\end{array}\right.
      \end{equation}

It is well-known that $\left\{w_k\right\}_{k\geq1}$ can be chosen
to constitute an orthonormal basis of $L^2((0, L); \rho)$ with the scalar product $(f, g)= \int_0^L f(x) g(x) \rho(x) dx$. The eigenvalues $\lambda_k$ are all simple, with $\lambda_1=0$ (and $w_1\equiv 1/[\int_0^L \rho (x) dx]^{1/2}$) and so that $\lambda_k$ grows quadratically with $k$. Similarly $|w_k(0)|\sim |w_k(0)|\sim 1$ uniformly on $k\ge 1$. This is a natural extension of the classical Sturm-Liouville problem for the 1-$d$ Laplacian in which case: $\lambda_k= (\pi (k-1))^2/L^2$ and $w_k(x)= \cos (\pi (k-1))x/L)$.

 On the basis of eigenfunctions of \eqref{eqdeltalambda}  solutions $\psi$ of \eqref{e45} can be expanded in Fourier series in the form
\begin{equation}\label{Fourier1}
 \psi(x, t) = \sum_{k \ge 1} \psi_k(t) w_k(x),
\end{equation}
 where
 $\psi_k(t)$ is the solution of the ODE
 $$
 \psi_k''(t)+ \lambda_k \psi_k(t)=a(L^-)w_k(L) s(t), \, 0< t <T; \quad  \psi_k(T)=  \psi_k'(T)=0,
 $$
 so that
 $$
 \psi_k(t)= -\frac{a(L^-)w_k(L)}{\sqrt{\lambda_k}} \int_{t}^T \sin (\sqrt{\lambda_k}(t-\tau)) s(\tau) d\tau.
 $$
 All in all,
  $$
 \psi(x, t) = -a(L^-) \sum_{k \ge 1} \frac{w_k(L)}{\sqrt{\lambda_k}}\int_{t}^T \sin (\sqrt{\lambda_k}(t-\tau)) s(\tau) d\tau\,  w_k(x),
 $$
 On the other hand, 
 \begin{equation}
\label{Fourier2b}
 \psi(0, t)= -a(L^-)\sum_{k \ge 1} \frac{w_k(0)w_k(L)}{\sqrt{\lambda_k}}\int_{t}^T \sin (\sqrt{\lambda_k}(t-\tau)) s(\tau) d\tau.
\end{equation}

 According to the observability inequality \eqref{e91} it follows that
\begin{equation}\label{e48b}
\left \Vert s(t)\right\Vert
_{H^{-1}\left(L\beta,T\right) } \leq C_1 \left \Vert a(L^-)\sum_{k \ge 1} \frac{w_k(0)w_k(L)}{\sqrt{\lambda_k}}\int_{t}^T \sin (\sqrt{\lambda_k}(t-\tau)) s(\tau) d\tau \right\Vert_{L_{2}\left(0,T\right) }
\end{equation}
for all $s \in H^{-1}\left(0,T\right)$ with support in $(L\beta,T)$.

Whether this estimate \eqref{e48b}, which we know is true by sidewise energy estimates, can be directly proved by Fourier series techniques is an interesting problem.
\medskip

 In the context of the control of 1-$d$ wave equation, often, the dual observability inequalities are proved using results in non-harmonic Fourier series whose ancestor is the classical Ingham inequality  mentioned above. 
This Ingham inequality has been the starting point of a rich literature motivated, to a large extent, by the dynamical and control theoretical properties of wave equations and to which Yves Meyer greatly contributed. 

Note however that the Ingham inequality, as such, can only be directly applied for the solutions of wave-like equations  considered here, but only when the  boundary conditions are homogeneous, i. e.
\begin{equation}\label{e45bxx}
\left\{~
\begin{aligned}
&\rho (x)p_{tt}-(a(x)p_{x})_{x}=0, && 0< x< L,~ 0< t<T\\
&p(x,T)=p_0(x), ~ p_{t}(x,T)=p_1(x),            && 0< x<L\\
&p_x(0,t)=0, ~ p_x(L,t)=0,                && 0< t<T
\end {aligned}
\right.
\end{equation} 
Solutions of \eqref{e45bxx}, when written in Fourier series in the basis $\{w_k\}$, take the form
\begin{equation}
 p(x, t) = \sum_{k \ge 1} \left [ p_{k, 0} \cos(\sqrt{\lambda_k}(T-t)) - \frac{p_{1,k}}{\sqrt{\lambda_k} }\sin(\sqrt{\lambda_k}(T-t) )\right ] w_k(x).
 \end{equation} 
 When considering their trace at $x=L$, they lead to expressions of the form  
 $$
\sum_{n\in\mathbb{Z}} a_n
e^{i\lambda_nt}. 
$$
 This allows to directly apply the Ingham inequality in this homogeneous case, since the time-frequencies $\{ \sqrt{\lambda_k} \}_{k\ge 1}$ fulfill the needed gap condition under the $BV$ assumptions on coefficients.

Whether non-homogeneous variants of the Ingham  inequalities can be developed so to achieve estimates of the form \eqref{e48b}, which we know are true thanks to the  sidewise energy estimates,
 is, as far as we know, an issue that has not been treated so far. 
 
 \section{Sidewise energy estimates by Fourier series}
 
One could try to use Fourier series expansions in a sidewise manner to recover the sidewise energy inequalities of the form \eqref{e48b}. When doing it  $x$ plays the role of time, while $t$ that of the space-variable. 
Accordingly, adapting the  notation, i. e. swapping the role of $x$ and $t$, so that  $x=\hat{t}$ and $t=\hat{x}$, equation \eqref{e45}  can then be written in the form
\begin{equation}\label{lat1}
(a(\hat{t}) q_{\hat{t}})_{\hat{t}}- \rho(t) q_{\hat{x}\hat{x}}=0,  \, 0<\hat{x}<T, 0<\hat{t}<L.
\end{equation}
This is a 1-$d$ wave equation in which coefficients depend solely on time $\hat{t}$ . Since $a$ and $\rho$ belong to $BV$, the equation, when complemented with suitable boundary conditions,  is well-posed in finite-energy spaces and enjoys the property of finite velocity of propagation.

On the other hand, dividing it by $\rho(\hat{t})$ it can readily be written in the form
\begin{equation}\label{lat2}
b(\hat{t}) q_{\hat{t}\hat{t}}- q_{\hat{x}\hat{x}} + c(\hat{t}) q_{\hat{t}}=0, \, 0<\hat{x}<T, 0<\hat{t}<L
\end{equation}
where $b(\hat{t})=a(\hat{t})/\rho(\hat{t})$ and $c(\hat{t})=a'(\hat{t})/\rho(\hat{t})$.

This equation, when considered in the full space, allows for a simply Fourier representation based on plane waves of the form
\begin{equation}\label{lat3}
q(\hat{x}, \hat{t})= q_{\xi}(\hat{t}) \exp(i\xi x)
\end{equation}
with $\xi \in \mathbb{R}$, where
$q_{\xi}= q_{\xi}(\hat{t}) $ fulfills
\begin{equation}\label{lat4}
b(\hat{t}) q_{\xi}''+ |\xi|^2 q_{\xi} +c(\hat{t}) q_{\xi}'=0.
\end{equation}

Energy estimates for \eqref{lat2} are equivalent to energy estimates for  frequency-dependent damped oscillators \eqref{lat2}.

In this setting, the adjoint system \eqref{e45} leads to a problem of the form \eqref{lat2} with initial, terminal,  and two boundary conditions of the form
\begin{equation}\label{lat5}
q_{\hat{t}}(\hat{x}, 0)=0,
\end{equation}
\begin{equation}\label{lat6}
q_{\hat{t}}(\hat{x}, L)=s(x),
\end{equation}
\begin{equation}\label{lat7}
q(T, \hat{t})=q_{\hat{x}}(T, \hat{t}) =0.
\end{equation}

This does not constitute a classical Cauchy problem. But, given that the source $s=s(\hat{x})$ entering in the condition \eqref{lat6}  has support in $(L\beta, T-L\beta)$, the system can be complemented with Dirichlet boundary conditions
\begin{equation}\label{lat8}
q(0, \hat{t})=q(T, \hat{t})=0,
\end{equation}
or also with the Neumann ones. 
Yet, to guarantee the solvability of the system, one would need to add a terminal condition at $\hat{t}=L$ of the form
$$
q(\hat{x}, L)=h(\hat{x})
$$
so to solve the equation backwards in the sense of time from $\hat{t}=L$ to $\hat{t}=0$.

Choosing $h$ with support in $(L\beta, T-L\beta)$, the finite velocity of propagation would guarantee that, other than \eqref{lat8}, the solution would also fulfill \eqref{lat7}.

But the added condition  \eqref{lat5} would impose an implicit compatibility condition on $h$.

Solutions of this auxiliary system could be represented in Fourier on the basis of eigenfunctions $\{\sin(k\pi \hat{x}/T)\}_{k \ge 1}$ and the observability inequality \eqref{e48b} would be a natural output of the well-posedness of this system in the positive sense of time, from  $\hat{t}=0$ to $\hat{t}=L$.

It would be interesting to see if this argument, which requires dealing with the implicit compatibility condition on $h$ to assure \eqref{lat5}, would eventually lead to an estimate of the form \eqref{e45}, as we expect. But, even if that were true, this would not be the approach inspired by Ingham inequalities we could aim for, as explained in the previous section, but rather a Fourier interpretation of the sidewise energy estimates above.

\section{Sidewise interior point-wise control}

We have shown that the sidewise controllability properties of the wave equation can be handled to show that, acting on one end, one can regulate the deformation or flux on the other one, Theorem \ref{controltheorem} and, simultaneously, also the final state, Theorem \ref{esc}.

The same questions can be formulated when one is trying to regulate the dynamics on an interior point $x_0\in (0, L)$ of the string. The problem would, for instance,   read as follows: {\it Given $z_0\in H^1(0, L)$ and $z_1\in L^2(0, L)$ and $p=p(t)\in H^1(0, T)$, to find a control $u=u(t)$  so that the solution of \eqref{e1} satisfies the terminal condition \eqref{final} and 
\begin{equation}\label{interiortraceb}
y(x_0, t)=p(t), 0\le t \le T.
\end{equation}
}
Of course, we could also consider the simpler problem in which we do not impose the final condition but simply \eqref{interiortraceb}.

As before, because of the finite velocity of propagation, one does not expect \eqref{interiortraceb} to hold when $T$ is too small. Indeed, the action of the control $u=u(t)$ needs a time of the order of $x_0\beta$ to get to the point $x=x_0$ where the dynamics is aimed to be controlled. But, in principle, in view of our previous results on sidewise control, one could expect this kind of result to hold when $T>2L\beta$. But, as we shall see, things get more complex when considering this point-wise interior control problem.

Using duality, once more, the problem can be transformed into the observability of the following adjoint system: 
\begin{equation}\label{e45b}
\left\{~
\begin{aligned}
&\rho (x)\psi_{tt}-(a(x)\psi_{x})_{x}=s(t)\delta_{x_0}, && 0< x< L,~ 0< t<T\\
&\psi(x,T)=\psi_0(x), ~ \psi_{t}(x,T)=\psi_1(x),            && 0< x<L\\
&\psi_x(0,t)=0, ~ \psi_x(L,t)=0,                && 0< t<T.
\end {aligned}
\right.
\end{equation} 

Note that in the present case the source $s=s(t)$ is located at $x=x_0$, and is represented by the Dirac delta $\delta_{x_0}$. And the goal this time is to observe the source $s=s(t)$ out of measurements done at $x=0$, together with the terminal conditions $(\psi_0, \psi_1)$, i.e. 
\begin{equation}\label{e48b'bb}
\left \Vert s(t)\right\Vert
_{H^{-1}\left(L\beta,T-L\beta\right) } + \left \Vert  \psi_0\right\Vert_{L^2(0, L)} + \left \Vert \psi_1 \right\Vert_{H^{-1}(0, L)}  \leq C \left \Vert \psi(0,t)\right\Vert
_{L^{2}\left(0,T\right) } 
\end{equation}

But as we shall see, \eqref{e48b'bb} may not hold.

\begin{remark} 
System \eqref{e45b} is the dual problem of the point-wise control problems considered in \cite{Lio3} and that, in the context of beam and plate equations, motivated the work of Haraux and Jaffard discussed in the introduction. There, the control was acting point-wise and the goal was to control the dynamics at the final time.

One of the first issues that arises in this context is the optimal regularity of solutions of equations of the form  \eqref{e45b} with a point-wise mass as source term, addressed in \cite{Lio3}, up to space dimensions $d \le 3$. The first proof of the optimal regularity in dimension $d=3$, which guarantees that, when $s=s(t)\in L^2(0, T)$, solutions lie in $C([0, T]; L^2_x)$ for the homogeneous Dirichlet boundary conditions on a bounded domain, is due to Yves Meyer \cite{Meyer3}, as reported in  \cite{Lio3} (see also \cite{Cohen}). The proof used H\"ormander - Pham estimates on the spectral function $e_\lambda(x)=\sum_{\lambda_n\le \lambda}|\varphi_n(x)|^2$, $\{\varphi_n\}_{n\ge 1}$ being the eigenfunctions of the Laplacian. Note that such a result is far from being achievable by classical energy estimates since the $C([0, T]; L^2_x)$ regularity of solutions of the wave equation requires that the right-hand side takes values in $L^1(0,T; H^{-1}_x)$. But $\delta_{x_0}$ does not belong to $H^{-1}_x$ beyond dimension  $d=1$.

\end{remark}

In the one-dimensional case under consideration, solutions of \eqref{e45b} are continuous across $x=x_0$, and they present a jump on the first derivative so that
$$
[a(x) \psi_x] (x_0, t) = a(x^+_0) \psi_{x}(x^+_0, t) - a(x^-_0) \psi_{x}(x^-_0, t)= -s(t).
$$

 One can apply sidewise energy estimates when making measurements both on $x=0$ and $x=L$. The following holds:
  \begin{prop}\label{counter}
Let us consider system \eqref{e45b} with coefficients satisfying the assumptions
\eqref{e43} and \eqref{e42}.

Let $T>2 M(x_0)\beta$ with $M(x_0)=\max (x_0, L-x_0)$. 

Then, there exists a constant $C>0$ such
\begin{equation}\label{e48b'bbb}
\left \Vert s(t)\right\Vert
_{H^{-1}\left(L\beta,T-L\beta\right) } + \left \Vert  \psi_0\right\Vert_{L^2(0, L)} + \left \Vert \psi_1 \right\Vert_{H^{-1}(0, L)}  \leq C \left (\left \Vert \psi(0,t)\right\Vert
_{L^{2}\left(0,T\right)} +  \left \Vert \psi(L,t)\right\Vert_{L^{2}\left(0,T\right) } \right)
\end{equation}
for every $s=s(t)$ such that $s\equiv 0$ away from $(M(x_0), T-M(x_0))$.
 \end{prop} 
 
But, as we shall see, sidewise observability from one side cannot be expected in this case. The counterexample reads as follows;
\begin{prop}
For all $x_0\in (0, L)$ and all $T>0$ there are  solutions of  \eqref{e45b} such that
\begin{equation}
\psi(0, t)=0, \quad 0<t<T; \quad s(t)\ne 0.
\end{equation}
In particular \eqref{e48b'bb} can never be achieved.
\end{prop}
\noindent {\bf Proof.} It is sufficient to define $\psi$ piecewise to the left and right of the point mass $x_0$:
\begin{itemize}
\item First, set $\psi\equiv 0$ in $(0, x_0) \times (0, T)$.
\item Second, take $\psi$ to be any non-trivial solution of 
\begin{equation}\label{e45bxx}
\left\{~
\begin{aligned}
&\rho (x)\psi_{tt}-(a(x)\psi_{x})_{x}=0, && x_0< x< L,~ 0< t<T\\
&\psi(x,T)=\psi_0(x), ~ \psi_{t}(x,T)=\psi_1(x),            && x_0< x<L\\
&\psi_x(x_0,t)=0, ~ \psi_x(L,t)=0,                && 0< t<T.
\end {aligned}
\right.
\end{equation} 
\item We then define $\psi$ piecewise, taking value $\psi\equiv 0$ on $0<x<x_0$ and $\psi$ as above in $x_0<x<L$. This generates a solution of \eqref{e45b} with a source
$s(t)= -a(x^+_0) \psi_x(x^+_0, t)$, which is non trivial  since $\psi\ne 0$ is a non-trivial solution of the wave equation on $ x_0< x<L$.
\end{itemize}
This concludes the proof of the Proposition.
\medskip

\begin{remark}
The presence of a Dirac mass at $x=x_0$ makes the systems under consideration to look like two interconnected strings, one at the left and the other one at the right of $x=x_0$, similar to those addressed in   \cite{LeugZ} and \cite{Dager}, in the context of waves in networks. But the spectral analysis in those works, which combined with irrationality and diophantine approximation tools, allows to get observability inequalities in weaker norms,  does not seem to apply here due to the counterexample Proposition \ref{counter}.
\end{remark}




\section{Multi-dimensional formulations}

As explained in \cite{Sarac}, all these problems make sense in the multi-dimensional context too, and not only for wave-like equations, but also for heat and Schr\"odinger ones, among others.

Controllability theory has been broadly developed in the multi-dimensional context. Duality arguments reduce the problem to boundary observability inequalities that can be obtained by different methods, including multipliers, Carleman inequalities, and microlocal analysis, \cite{survey}.
However, the adaptation of the existing techniques to handle sidewise problems seems to be far from obvious.

Let us present this interesting and challenging issue in the context of the multi-dimensional wave equation. We refer the interested readers to \cite{DZ} for a systematic analysis of this problem using microlocal tools.

Let $\Omega$ be a bounded smooth domain of $\mathbb{R}^d$, in dimension $d\ge 2$, and consider the wave equation:
\begin{equation}\label{e1multi}
\left\{~
\begin{aligned}
&y_{tt}-\Delta y =0,  &&\hbox{ in } \Omega \times (0, T)\\
&y(x,0)=y_0(x), ~y_{t}(x,0)=y_1(x),&& \hbox{ in } \Omega\\ 
&\partial y/\partial \nu=u,              &&\hbox{ on } \Gamma_c\times (0, T) \\
&\partial y/\partial \nu=0,                 &&\hbox{ on } \Gamma_0\times (0, T).
\end {aligned}
\right.
\end{equation} 

Here $\Gamma_0$ and $\Gamma_c$ stand for a partition of the boundary, $\Gamma_0$ being the  part of the boundary without any applied force, and $\Gamma_c$ the one under control. Here and in what follows $\nu$ denotes the outward unit normal vector and $\partial \cdot/\partial \nu$ the normal derivative.

Given a smooth enough target $p: \Gamma_0\times(0, T)\to \mathbb{R}$  to be tracked, the question is then to find a control $u$ in, say, $L^2(\Gamma_c \times (0, T))$, such that the solution $y$ of \eqref{e1multi} satisfies the condition
$$
y = p \quad \hbox{on} \ \Gamma_0\times(0, T).
$$

This is the natural multi-d version of the sidewise controllability problem discussed above.

It is equivalent, by duality, to a new class of sidewise observability inequalities for the adjoint system
\begin{equation}\label{e1multi2}
\left\{~
\begin{aligned}
&\psi_{tt}-\Delta \psi =0,  &&\hbox{ in } \Omega \times (0, T)\\
&\psi(x,T)=\psi_{t}(x,T)=0,&& \hbox{ in } \Omega\\
&\partial \psi/\partial \nu=0,              &&\hbox{ on } \Gamma_c\times (0, T) \\
&\partial \psi/\partial \nu=s,                 &&\hbox{ on } \Gamma_0\times (0, T).
\end {aligned}
\right.
\end{equation} 
Here $s=s(x, t)$ is a smooth boundary condition given on $\Gamma_0\times (0, T)$. The question is then  whether one can prove the existence of an observability constant $C>0$ such that
\begin{equation}
||s||^2_{*} \le C ||\psi||_{L^2(\Gamma_c\times (0, T))}^2
\end{equation}
for every solution of this adjoint system.

Observe that, here, the norm $||\cdot ||_*$ in this inequality is to be identified both in what corresponds to the Sobolev regularity and the support within $\Gamma_0 \times (0, T)$.

Techniques based on non-harmonic Fourier series and Ingham-type inequalities cannot be employed in the multi-dimensional context. This is simply due to the fact that, according to Weyl's Theorem, eigenvalues of the laplacian and more general second elliptic operators grow, roughly,  as $\lambda_k \sim C(\Omega) k^{2/d}$ when $k \to \infty$. Accordingly, $\sqrt{\lambda_k}\sim C(\Omega) k^{1/d}$ cannot fulfill the needed gap condition or the more general asymptotic density properties of the spectrum that the application of Beurling-Malliavin's theorem requires, \cite{HJ}.

 As we mentioned above, the existing techniques do not seem to yield this kind of multi-d sidewise observability inequalities in a direct manner.
 
 However, as described in \cite{L2}, using Holmgren's uniqueness Theorem, a unique continuation property can be easily proved. This constitutes a weaker and non-quantitative version of this kind of inequality. 
 
Indeed, by Holmgren's Theorem, as soon as $ \psi \equiv 0$ in $\Gamma_c\times (0, T)$ and $T$ is large enough, one can guarantee that $s\equiv 0$ provided its support is localized in a subset of the boundary of the form $\gamma \times (\tau, T-\tau)$, with $\gamma$ a suitable open subset of $\Gamma_0$ and $0 < \tau <T/2$. Essentially, $\gamma$ is constituted by the points for which the geodesic distance (within $\Omega$)  to $\Gamma_c$ is less than $\tau$. 
 
 Obviously, for this result to be active in some effective subset $\gamma \times (\tau, T-\tau)$, one needs $T$ to be large enough, in particular $T>2\delta$, where $\delta$ is the minimal geodesic distance from $\Gamma_c$ to $\Gamma_0$ (\cite{Cazenave}, \cite{Haraux}).

 This unique continuation result assures, in the corresponding geometric setting, that the wave equation enjoys the property of sidewise approximate controllability: i. e. that given any $p\in L^2(\gamma \times (\tau, T-\tau))$ and any $\epsilon >0$ there exists a control $u\in L^2(\Gamma_c \times (0, T))$ (depending on $\epsilon$) such that the corresponding solution $y$ satisfies
$$
||y  - p ||_{L^2(\gamma \times (\tau, T-\tau))} \le \epsilon.
$$

This result can be viewed as a partial extension of the 1-$d$ results in this paper  to the multi-dimensional case. However, these arguments, based purely on Holmgren's uniqueness Theorem, do not yield any quantitative estimates.

Of course, this discussion can also be extended to the problem  in which we aim at both simultaneously: using the control of the trace $y$ over $\gamma$ to track a given function $p$, but also the final value of the solution at $t=T$.

The systematic analysis of these problems in the multi-dimensional context for the wave equation and other relevant models constitutes a very rich source of interesting open problems. We refer to \cite{DZ} for an in-depth discussion employing microlocal tools. Whether the theory of non-harmonic Fourier series theory can contribute to treating these problems for platelike models as in \cite{HJ} constitutes an interesting issue worth to be investigated.

\bigskip

\section*{Acknowledgments} The author thanks A. Bonami, A. Haraux, E. Tr\'elat and the anonymous referees for their useful remarks.

The author has been funded by the Alexander von Humboldt-Professorship program, the ModConFlex Marie Curie Action, HORIZON-MSCA-2021-$d$N-01, the COST Action MAT-DYN-NET, the Transregio 154 Project ``Mathematical Modelling, Simulation and Optimization Using the Example of Gas Networks'' of the DFG, grants PID2020-112617GB-C22 and TED2021-131390B-I00 of MINECO (Spain), and by the Madrid Goverment -- UAM Agreement for the Excellence of the University Research Staff in the context of the V PRICIT (Regional Programme of Research and Technological Innovation).

\bibliographystyle{amsalpha}

\end{document}